\documentclass{amsart}
\usepackage{amsmath}
\usepackage{amssymb}
\usepackage{amsopn}
\usepackage{amscd}
\usepackage{epsfig}
\usepackage{amsfonts}
\usepackage{latexsym}
\usepackage{graphicx}

\usepackage{pdfsync}

\usepackage{tikz}
\usetikzlibrary{positioning,arrows}

\newcommand{\ca}{c}
\newcommand{\dd}{\check{1}}
\newcommand{\ee}{\check{2}}

\newcommand{\thetaSF}{\theta^{^{\scriptscriptstyle  \diamond}}}

\theoremstyle{plain}
\newtheorem{theorem}{Theorem}[section]

\newtheorem{lemma}{Lemma}[section]
\newtheorem{proposition}{Proposition}[section]
\newtheorem{procedure}{Procedure}[section]
\newtheorem{corollary}{Corollary}[section]

\newtheorem{problem}{Problem}[section]

\title[]{Topological conjugacy of constant length substitution dynamical systems}

\author[]{Ethan M.~Coven}
\author[]{F.~Michel Dekking}
\author[]{Michael S.~Keane}

\address{Wesleyan University and Delft University of Technology.}

\date{\today}

\begin{document}


\maketitle

\begin{abstract}
 Primitive constant length substitutions generate minimal
symbolic dynamical systems. In this article we
present an algorithm which can produce the list of injective substitutions
of the same length that generate topologically conjugate systems.
 We show that each conjugacy class contains infinitely substitutions  which are not injective.
As examples, the Toeplitz conjugacy class contains three injective
substitutions (two on two symbols and one on three symbols), and the
length two Thue-Morse conjugacy class contains twelve
substitutions, among which are two on six symbols. Together, they
constitute a list of all primitive substitutions of length two with
infinite minimal systems which are factors of the Thue-Morse system.

\smallskip


\medskip


\end{abstract}

\section{Prologue}

In the article [3] published in 1971, the minimal dynamical systems arising
from primitive substitutions on a binary alphabet having the same constant
length were classified, yielding for a given such substitution a list of all
substitutions of the same length generating topologically conjugate systems.
Here we extend this classification to arbitrary finite alphabets.
More recently, the articles [4] and [5] exhibit
characterizations of such systems; these only implicitly yield corresponding
topological conjugacies, and do not result in lists of conjugate systems. Also, in \cite{HP} and \cite{Mentz} a related goal has been partially accomplished
---a classification of measure-theoretic conjugacy---for a restricted class of constant length substitutions.

If two constant length substitution systems are topologically conjugate, then the lengths of the substitutions are powers of the same integer (\cite{Dur2008},\cite{CDL2014}). Therefore,  by taking suitable powers we can, and do,  restrict our attention to substitutions of the same length $L$.

 In this contribution  we address  the following two problems, in which $L$ denotes a fixed integer larger than one.

\begin{problem}\label{prob:a}
 Let $\alpha$ and $\beta$ be two substitutions of the same length $L$, both primitive.
Decide whether the dynamical systems $(X_\alpha,\sigma)$ and $(X_\beta,\sigma)$ are topologically conjugate.
\end{problem}

\begin{problem}\label{prob:b}
 Let $\alpha$  be a primitive substitution of length $L$.
Give a list of all the injective substitutions $\beta$ of  length $L$ such that the dynamical systems
$(X_\alpha,\sigma)$ and $(X_\beta,\sigma)$ are topologically conjugate.
\end{problem}

Finite systems are elementary, and we restrict attention everywhere
to the non-periodic case of primitive substitutions with corresponding infinite minimal sets.

 We show  that to any primitive
substitution of constant length whose minimal set is infinite, there are always
infinitely many primitive substitutions of the same constant length having
topologically conjugate minimal systems, but only finitely many of these are
injective. Thus, the list produced by our algorithm for attacking Problem 1.2 will,
 starting from any given primitive substitution of constant length, consist of all
injective substitutions of that length with dynamical systems topologically conjugate
to the initial system.
Clearly, since the list in Problem~\ref{prob:b} is finite, Problem~\ref{prob:a} has then also
been solved, since there is a simple algorithm to associate to a substitution an injective substitution generating a conjugate system
(cf. Section \ref{sec: injective}). This contrasts with the situation for the natural generalization  of our problem to the
collection of \emph{all} substitutions. In \cite{Dek2014} it is shown that there may be infinitely many
primitive injective (non-constant length) substitutions that generate systems conjugate to a
system generated by a substitution with the same Perron-Frobenius eigenvalue for its incidence matrix.


 Recently  a completely different solution has been obtained  for  Problem \ref{prob:a}. in the preprint \cite{CQY2015}. Actually, because of Theorem \ref{th: three-block}, a solution of Problem \ref{prob:a} also yields a solution of Problem \ref{prob:b}. However, it seems unfeasible---using the algorithm of \cite{CQY2015}---to obtain the Thue-Morse list by hand, as we do in Section \ref{sec: morse}.

\section{Substitutions and standard forms}\label{sec: standard}

We begin by recalling the basic definitions and known results without proof for primitive substitutions
and their corresponding minimal systems, referring the reader to the standard reference \cite{Que}.

\medskip

Let $A$ be a finite set (an {\it alphabet}) with $c\geq2$ elements which are {\it symbols}, or {\it letters}.
Elements of $A^*=\cup_{n=0}^\infty A^n$ are called {\it words}. A {\it substitution} is a mapping
$$
\alpha:A\longrightarrow A^*.
$$
The substitution $\alpha$ is of constant length $L$ if $\alpha (a)\in A^L$  for each $a\in A$. It is natural to
view $A^*$ as a semigroup under juxtaposition, thus extending $\alpha$ to mappings from $A^*$ to $A^*$, $A^{\bf N}$
to $A^{\bf N}$, and $A^{\mathbb{Z}}$ to $A^{\mathbb{Z}}$ - no confusion results if we also denote them by $\alpha$, and they
can be iterated, defining $\alpha^n$ for each $n\in \bf N$.

\smallskip

\noindent {\bf Definition.} The substitution $\alpha$ is {\it primitive} if for some  $n>0$ and for every $a\in A$
the word $\alpha^n(a)$ contains each of the letters of $A$. The {\it language} of $\alpha$ is
the subset ${\mathcal L}_\alpha$ of $A^*$ consisting of those words appearing as consecutive letters, {\it subwords}, or {\it factors},
of images under powers of $\alpha$.
We denote by ${\mathcal L}_\alpha^N$ the set of  words of length $N$ in ${\mathcal L}_\alpha$.

\smallskip

We write $X_\alpha$ for the compact subset of $A^{\mathbb{Z}}$ of bilaterally infinite sequences each of whose finite factors
belongs to the language of $\alpha$. Under the left shift $\sigma$ on $A^{\mathbb{Z}}$, it is a minimal symbolic system
whenever $\alpha$ is primitive. If in addition, $X_\alpha$ is infinite, then $\alpha$ is {\it recognizable} (\cite{Mos}).
For constant length $L$ substitutions, this is equivalent to the existence of a  semi-conjugacy from the minimal
system $(X_\alpha,\sigma)$ to the rotation by $1$ on the compact group of $L$-adic integers, which describes a
unique hierarchical structure for each of the sequences belonging to $X_\alpha$.

\smallskip

For  substitutions, it is clear that the names we give to the individual symbols of their
alphabets are not essential - different namings will produce conjugate systems. This leads us to restricting
an alphabet of $\ca$ symbols to the alphabet $A=\{1,\dots,\ca\}$. Even then, there is a permutational ambiguity,
since permuting $A$ will yield up to $\ca!$ different substitutions, which we view as essentially the same. We
find it useful in the following to single out one of these permutations as the one yielding the {\it standard
form} of a substitution, as follows. If $\alpha$ is a constant length $L$ substitution on the alphabet of size
$\ca$, then we define its {\it characteristic word} to be the word $\alpha(1)\cdots \alpha(\ca)$ of length $L\ca$.
For constant length\footnote{See \cite{Dek2015} for a standard form for arbitrary substitutions.} substitutions, permutations yielding different substitutions then possess different characteristic words, and we call the
substitution with the lexicographically smallest characteristic word the {\it standard form} of the substitution
$\alpha$.

\section{Letter-to-letter maps}\label{sec:L2L}

Let $A$ and $B$ be finite alphabets. A map
$$
\pi:A\longrightarrow B
$$
is called a {\it letter-to-letter map\,}; by juxtaposition it clearly extends to maps from (finite
or infinite) sequences on $A$ to sequences of the same lengths on $B$. We also denote this
extension by the same symbol $\pi$. It will appear that the following easily proved lemma is the key to
understanding the properties of conjugacies.

\begin{lemma}\label{lem: twining}  If $\alpha:A\rightarrow A^*$ and $\beta:B\rightarrow B^*$ are substitutions, and
if $\pi$ satisfies the { \it intertwining} equation $\pi\,\alpha=\beta\,\pi$, then for each positive integer $n$
$$
\pi\,\alpha^n=\beta^n\,\pi.
$$
\end{lemma}

%

Under the hypotheses of the lemma, the word $\alpha^n(a)$ is mapped by $\pi$ to the word $\beta^n(b)$, with
$b=\pi(a)$, for any positive $n$. In particular, the language of $\alpha$ is mapped to the
language of $\beta$, and we have:

\vskip.5cm

\begin{corollary}\label{corr:min}
 $\pi(X_\alpha)\subseteq X_\beta$, with equality whenever $\pi$ is surjective.
In particular, if $\pi$ is surjective, then primitivity of $\alpha$ implies primitivity of
$\beta$ and minimality of $(X_\alpha,\sigma)$ implies minimality of $(X_\beta,\sigma)$.
\end{corollary}

 When $\pi\,\alpha=\beta\,\pi$ and $\pi$ is surjective, we call $\beta$ an \emph{amalgamation} of $\alpha$.

\section{$N$-Block presentations and $N$-Block substitutions}\label{sec:Hedlund}


\vskip.2cm
Let $A$ be a finite alphabet, and let $N\ge 2$ denote a positive integer.
We consider the elements $a_0a_1\dots a_{N-1}$ of $ A^N$ as symbols in an alphabet denoted $A^{[N]}$ by  defining the $N$-block map
$$\Psi(a_0a_1\dots a_{N-1})=[a_0a_1\dots a_{N-1}].$$
If $X$ is a closed $\sigma$-invariant subset of $A^{\mathbb{Z}}$, then $X^{[N]}:=\psi(X)$ is called the {\it $N$-block presentation} of $X$, where $\psi$
is the conjugacy from $(X,\sigma)$ to $(X^{[N]},\sigma)$ associated to the sliding block code $\Psi$ (see e.g.~\cite{Will}).
The inverse of $\psi$ is associated to the letter-to-letter map $\pi_0$ given by
$$\pi_0([a_0a_1\dots a_{N-1}])=a_0.$$

\medskip

We now concentrate our attention on $X=X_\alpha$, where $\alpha$ is a  primitive substitution on $A$ with constant length $L$.
 A pleasant property is that the $N$-block presentation of $(X_\alpha,\sigma)$ is again a substitution dynamical system. If we define the alphabet $B=A^{[N]}_\alpha:=\{[a_0\dots a_{N-1}]: a_0\dots a_{N-1}\in \mathcal{L}_\alpha$\}, then the map $\pi_0$ from $B$ to $A$ satisfies
the intertwining condition $\pi_0\beta=\alpha\pi_0$ of the previous section, if we define the substitution
$\beta$ on $B$ properly. Moreover, since $\pi_0$ is obviously surjective, Corollary \ref{corr:min} then implies that
the systems $(X_\alpha, \sigma)$ and $(X_\beta, \sigma)$ are conjugate.

Such a $\beta$ exists, and has been introduced in  Queff\'elec's book \cite{Que} on page 95. However, we want a whole family of substitutions generating the $N$-block presentation of $(X_\alpha,\sigma)$. We denote the members of this family by $\widehat{\alpha}_{N,M}$. Here $M$ is called the {\it lag} of $\widehat{\alpha}_{N,M}$. The substitutions $\widehat{\alpha}_{N,0}$ are considered in \cite{Que}, and the $\widehat{\alpha}_{2,M}$ play a key role in \cite{HP}.

If $[a_0\dots a_{N-1}]$ is an element of $B$, we can  apply $\alpha$, obtaining a word
$$v=v_0v_1\dots v_{LN-1}:=\alpha(a_0\dots a_{N-1}).$$
Now choose any integer $M$ with $0\leq M \leq (L-1)(N-1)$, so that the factor $w$ of length $L+N$  of $v=\alpha(a_0\dots a_{N-1})$
starting with the symbol $v_M$ is well--defined. Then we define
$$\widehat{\alpha}_{N,M}([a_0\dots a_{N-1}])= [v_M\dots v_{M+N-1}][v_{M+1}\dots v_{M+N}]\dots[v_{M+L-1}\dots v_{M+L+N-2}].$$

\smallskip

\noindent {\bf Example.} Let $A=\{1,2,3\}$, and let $\alpha$ be given by
$$\alpha(1)=1233,\; \alpha(2)=2313,\;\alpha(3)=3123.$$
Then the words of length $N=2$ in the language of $\alpha$ are $12, 13, 23, 31, 32$ and $ 33$.  We construct the 2-block substitution $\beta=\widehat{\alpha}_{2,1}$ on the alphabet $A^{[2]}_\alpha$ with lag $M=1$.
Since $\alpha(12)=12332313$, we have $\beta([12])=[23][33][32][23]$. Coding the $[aa']$ in lexicographical order to a standard alphabet gives $B=\{1,2,\dots,6\}$. On\footnote{In the sequel we will often identify the alphabet $A^{[N]}_\alpha$ with its standard form.} $B$ we have $\beta(1)=3653,\beta(2)=3664,\; \beta(3)=4264,\;\beta(4)=1341,\; \beta(5)=1353,\;\beta(6)=1364.$

\smallskip

\begin{proposition} Let $\alpha$ be a primitive substitution of length $L$ on an alphabet $A$. For a positive  integer $N$, and any $M$ with $0\leq M \leq (L-1)(N-1)$ let $\beta=\widehat{\alpha}_{N,M}$ on the alphabet $A^{[N]}_\alpha$.
Then the system $(X_\beta,\sigma)$  is conjugate to the system $(X_\alpha,\sigma)$.\label{th:thetahat}
\end{proposition}

{\it Proof:} (From \cite{HP}.) One can show that $\beta=\widehat{\alpha}_{N,M}$ on $B=A^{[N]}_\alpha$ is a primitive substitution.
 One easily verifies that under the projection $\pi_0$ we have $\pi_0(X_\beta)\subseteq X_\alpha$. Then by minimality, the sets are equal. \qed

\smallskip

An alternative proof for this proposition can be given using the following lemma. For notational reasons we define the {\em hat operator} $\mathcal{H}$ by
$\mathcal{H}_{N,M}(\alpha)=\widehat{\alpha}_{N,M}.$

\begin{lemma}
For all $n\ge 1$, $[\mathcal{H}_{N,M}(\alpha)]^n=\mathcal{H}_{N,M(L^n-1)/(L-1)}(\alpha^n).$
\end{lemma}

{\it Proof:} It is easily seen that for two lags $M$ and $M'$ we obtain for the composition $\mathcal{H}_{N,M}(\alpha)\circ \mathcal{H}_{N,M'}(\alpha)=\mathcal{H}_{N,M'L+M}(\alpha)$.
Iterating $\alpha$, the cumulative lag in $\alpha^n(a_0\dots a_{N-1})$ is $L^nM+L^{n-1}M+\dots+LM+M=M(L^n-1)/(L-1)$.\qed

\smallskip

If $X_\alpha$ is infinite, then clearly the alphabets $A^{[N]}_\alpha$ grow larger and larger with $N$. So by Proposition \ref{th:thetahat} one obtains

\begin{theorem}  For any primitive constant length substitution with infinite associated symbolic system there
exist infinitely many  primitive substitutions of the same length with symbolic systems
topologically conjugate to the given system.
\end{theorem}


\section{For substitution minimal sets 3-block codes suffice}\label{sec: 3-block}

 In general a semi-conjugacy from a system $(X,\sigma)$ to $(Y,\sigma)$ can always be obtained as a  sliding block code from $X$ to $Y$
 (see \cite{Hed}).

Here we give a new proof of a known result (see \cite{CDKL2014}, Theorem 3).

\smallskip

\begin{theorem}\label{th: three-block}
Let $\alpha$ and $\beta$ each be primitive  substitutions
of constant length $L>1$, whose minimal systems $(X_\alpha,\sigma)$ and
$(X_\beta,\sigma)$ are infinite.
 If there exists a semi-conjugacy from $(X_\alpha,\sigma)$ to
$(X_\beta, \sigma)$, and $\beta$ is injective then there is such a semi-conjugacy which is given by a 3-block   code.
\end{theorem}

\emph{Proof:} Denote by $\phi$ the hypothesized semi-conjugacy.
 We may assume without loss of generality that the associated sliding block code $\Phi$ is an $L^n$-block code with memory 0 for some integer $n$.

Recall that $\mathcal{L}^3_\alpha$ denotes the set of  words of length three in ${\mathcal L}_\alpha$, and let $B$ be the alphabet of $\beta$.
The proof now consists of two steps:
\vskip.2cm
Step 1. Construction of a three-block code $\Psi$ from $\mathcal{L}^3_\alpha$ to $B$.

Choose any three-block $ijk \in \mathcal{L}^3_\alpha$. The block $\alpha^{n}(ijk)$
is a $3L^n$-block from the language of $\alpha$, to which we can apply $\Phi$,
obtaining a $(2L^n\!+\!1)$-block of $\beta$ in ${\mathcal L}_\beta$.

By recognizability, there is a unique $\beta^n$-block, say $\beta^n(p)$,
occurring at a fixed position (independent of the choice of $ijk$) in this
block.
By injectivity of $\beta^n$,  $\Psi(ijk):= p$ is then well-defined.
\vskip.2cm
Step 2. The  block code $\Psi$ defines a map $\psi$ from $X_\alpha$ to a
closed, shift-invariant set $Y$ of sequences from the alphabet $B$, so that $\psi$
is a semi-conjugacy from $(X_\alpha,\sigma)$ to $(Y,\sigma)$. We show in this step that
$Y=X_\beta$.

To verify this, choose any $x\in X_\alpha$, apply $\alpha^n$ to $x$, then apply $\phi$,
and finally ``decode" using recognizability of $\beta^n$. The resulting sequence must then
be an element of $X_\beta$, and by minimality all elements of this set occur.
\qed

\medskip

{\bf Corollary.} If the semi-conjugacy of the  3-block Theorem is a conjugacy, then the 3-block
code which results from the proof is also a conjugacy.

\vskip.1cm

\emph{Proof:} If $x$ and $x'$ are different points in $X_\alpha$, it is obvious that their
images under $\psi$ are also different, so that a conjugacy results.
\qed

\smallskip

{\bf Remark.} In \cite{HP} it is shown for a rather special class of substitutions that the measure-theoretic semi-conjugacies are given by 2-block codes.
The example of the Thue-Morse substitution (see Section \ref{sec: morsefac}) shows that 3-block codes are sometimes necessary.

\section{Injective substitutions}\label{sec: injective}

A key ingredient in our classification result is that we may suppose that the substitutions are injective.
This is based on the following result.

\begin{theorem} {\rm ( \cite{BDM2004})}\label{th: injective}
Any system generated by a primitive, non-periodic substitution which is not injective is conjugate to a system generated by a primitive substitution that \emph{is} injective.
\end{theorem}

The proof given in \cite{BDM2004} is constructive, and yields what we call the {\it  injectivization} of a substitution.
It is an amalgamation of the original substitution.
The construction amounts to identifying (iteratively) those letters which have equal images.
For example, the substitution $\beta$ given by $$\beta(1)=46,\; \beta(2)=45,\;\beta(3)=26,\; \beta(4)=25,\;\beta(5)=13,\; \beta(6)=13$$  amalgamates in a first step to
$$\beta'(1)=45,\; \beta'(2)=45,\;\beta'(3)=25,\; \beta'(4)=25,\;\beta'(5)=13,$$
and then in a second step to the injective substitution
$$\beta''(1)=35,\; \beta''(3)=15,\; \beta''(5)=13.$$


\section{Substitutions and graph homomorphisms}\label{sec: graph}

Let $x$ be an infinite two-sided sequence over an alphabet $A$. Here we study the general question whether $x$ can be generated by a substitution of length $L$.

We  consider graphs ${\mathcal G}=(V,E),\; {\mathcal G'}=(V',E')$, and graph homomorphisms $\varphi:{\mathcal G}\rightarrow {\mathcal G'}$, i.e., maps
$\varphi: V\rightarrow V'$ having the property that $(u,v)\in E$ implies that $(\varphi(u),\varphi(v))\in E'$.

Let $W_2=W_2(x)=\{ab: ab=x_k x_{k+1}  \;\mathrm{for\, some}\; k \in \mathbb{Z}\}$, be the set of 2-blocks occurring in $x$, and  for  $0\le M \le L-1$ let
$W_{L,M}=W_{L,M}(x)=\{a_1\dots a_L: a_1\dots a_L=x_{kL+M}\dots x_{kL+M+L-1} \;\mathrm{for\, some}\; k \in \mathbb{Z} \}$ be the set of of $L$-blocks occurring in $x$ at positions $M \mod L$.

With $x$ we associate a family of graphs--- cf. \cite{Lothaire}, Section 1.3.4. The simplest is ${\mathcal G}_1^x=(V_1,E_1)$, the factor graph of order 1 of $x$, given by
$$ V_1=A, \quad E_1=\{(a,b): ab\in W_2\}.$$
The graphs ${\mathcal G}_{L,M}^x=(V_{L,M},E_{L,M})$  for $ M=0,\dots, L-1$ are defined by
$$V_{L,M}=W_{L,M},\quad E_{L,M}=\{(a_1\dots a_L,b_1\dots b_L): a_1\dots a_Lb_1\dots b_L\in W_{2L,M} \}.$$
We follow the convention of calling a surjective homomorphism an {\it epimorphism}. This requires that both the map on vertices and the map on edges are surjective.

\begin{lemma}\label{th: graph}
Let $x$ be sequence over $A$, and let $\varphi$ be a primitive substitution of length $L$ over $A$.  If $x$ is in $X_\varphi$  then
$\varphi$ is a  graph epimorphism, $\varphi: {\mathcal G}_1^x\rightarrow {\mathcal G}_{L,M}^x$ for some  $0\le M \le L-1$.
\end{lemma}

\emph{Proof:}  When $x$ is in $X_\varphi$,  $x$ can be written as a concatenation of  $\varphi$-blocks. Define $M$ as the first cutting position at or after 0.  Let $y$ be such that $x=\sigma^M\varphi(y)$. By minimality of $X_\varphi$, all letters of $A$ occur in $y$, and the substitution defines a surjective map from $V_1=A$ to $V_{L,M}=W_{L,M}(x)$. By minimality of $X_\varphi$,
 one has $W_2(y)=W_2(x)$, and if $ab$ occurs in $y$, then $\varphi(a)\varphi(b)$ is in $W_{2L,M}(x)$.   Thus  $\varphi$  can be seen as a graph homomorphism, and is also surjective on the edges, since any $w\in W_{2L,M}(x)$  must come from (at least) one word $ab$ in $W_2(y)$ as $w=\varphi(ab)$. \qed

\smallskip

Note that to avoid cumbersome notation we do not distinguish between $\varphi$ as a map on words and $\varphi$ as a graph epimorphism.


\smallskip

\noindent{\bf{Example: Thue-Morse sequence}}.

We consider the Thue-Morse sequence $x=0110100110010110\dots$.
It is easy to write down the graphs of the letters and the 2-blocks:


\smallskip

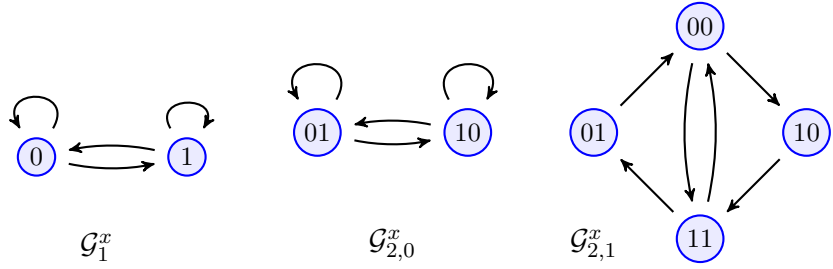
\begin{tikzpicture}[scale=.8,->,>=stealth',shorten <=8pt, shorten >=8pt,auto,node distance=2cm, thick,
  main node/.style={circle,fill=blue!10,draw,circle,inner sep=0pt},
  loop/.style={min distance=10mm,looseness=15,in=120,out=60},
  loop clock/.style={min distance=10mm,looseness=10,in=60,out=120}]
  \node[main node] (na) {$0$};
  \node[main node] [right of=na] (nb) {$1$};
  \draw[blue,thick,fill=blue!8] (na) circle [radius=.31cm] (nb) circle [radius=.31cm];
  \node at (na)  {$0$}; \node at (nb)  {$1$};
  \draw[->]
    (na) edge [bend right=12] (nb) (nb) edge [bend right=12] (na)
    (na) edge  [loop] node {} (na) (nb) edge  [loop clock] node {} (nb);
  \node  at (1,-1.5) {{\large ${\mathcal G}_1^x$}};
\end{tikzpicture}\label{G1-Morse}\hspace*{.5cm}
\begin{tikzpicture}[->,>=stealth',shorten <=8pt, shorten >=8pt,auto,node distance=2cm, thick,
  main node/.style={circle,fill=blue!10,draw,circle,inner sep=0pt},
  loop/.style={min distance=10mm,looseness=10,in=120,out=60},
  loop clock/.style={min distance=10mm,looseness=10,in=60,out=120}]
  \node[main node] (na) {$01$};
  \node[main node] [right of=na] (nb) {$10$};
  \draw[blue,thick,fill=blue!8] (na) circle [radius=.31cm] (nb) circle [radius=.31cm];
  \node at (na)  {$01$}; \node at (nb)  {$10$};
  \draw[->]
    (na) edge [bend right=12] (nb) (nb) edge [bend right=12] (na)
    (na) edge  [loop] node {} (na) (nb) edge  [loop clock] node {} (nb);
  \node  at (1,-1.5) {{\large ${\mathcal G}_{2,0}^x$}};
\end{tikzpicture}\label{G20-Morse}\hspace*{.5cm}
\begin{tikzpicture}[->,>=stealth',shorten <=8pt, shorten >=8pt,auto,node distance=2cm, thick,
  main node/.style={circle,fill=blue!10,draw,circle,inner sep=0pt},
  loop/.style={min distance=10mm,looseness=10,in=120,out=60},
  loop clock/.style={min distance=10mm,looseness=10,in=60,out=120}]
  \node[main node] (na) {$01$};
  \node[main node] [above right of=na] (nb) {$00$};
  \node[main node] [below right of=nb] (nd) {$10$};
  \node[main node] [below right of=na] (nc) {$11$};
  \draw[blue,thick,fill=blue!8] (na) circle [radius=.31cm] (nb) circle [radius=.31cm] (nc) circle [radius=.31cm] (nd) circle [radius=.31cm];
  \node at (na)  {$01$}; \node at (nb)  {$00$};  \node at (nc)  {$11$}; \node at (nd)  {$10$};
  \draw[->]
    (nb) edge [bend right=12] (nc) (nc) edge [bend right=12] (nb)
    (na) edge (nb) (nb) edge (nd) (nd) edge (nc) (nc) edge (na) ;
  \node  at (0,-1.5) {{\large ${\mathcal G}_{2,1}^x$}};
\end{tikzpicture}\label{G21-Morse}

\medskip

Note that ${\mathcal G}_{2,1}^x$ has too many vertices, and with ${\mathcal G}_{2,0}^x$ we find \emph{two} surjective graph homomorphisms: $\varphi(0)=01,\,\varphi(1)=10$,
corresponding to the usual substitution, but also $\varphi^\flat$ given by $\varphi^\flat(0)=10,\,\varphi^\flat(1)=01$. Note that both are in standard form.

%
%

\section{The list problem}\label{sec: sol}

  In this section we first describe an algorithm to find  for a given primitive substitution $\alpha$  all primitive injective substitutions $\beta$ of the same length whose associated systems are factors of $(X_\alpha,\sigma)$.

\begin{procedure}\label{th: factor}

   By Theorem \ref{th: three-block} we may suppose that the factor map is a 3-block map.
 Start with the 3-block presentation $X_{\alpha}^{[3]}$ of $\alpha$
from Section \ref{sec:Hedlund}. All factors of $(X_\alpha,\sigma)$ can be obtained by going through  all (including the identity)  letter-to-letter maps $\pi$ from $X_{\alpha}^{[3]}$  to another shift space. To see whether such a factor $X:=\pi(X_{\alpha}^{[3]})$ is generated by a primitive substitution of length $L$,  take any sequence $u$ from $X_{\alpha}^{[3]}$, and define $x:=\pi(u)$.
 Determine the graph ${\mathcal G}_1^x$ and the graphs ${\mathcal G}_{L,M}^x$ for all $M=0,...,L-1$. Then determine all epimorphisms $\varphi$ from  ${\mathcal G}_1^x$ to ${\mathcal G}_{L,M}^x$. By Lemma~\ref{th: graph} this gives a list of all possible candidates $\varphi$ that might generate $X$. Discard the $\varphi$ which are not primitive. Then check whether all subwords that appear in sequences of $X$ also occur in  sequences of $X_\varphi$.
 If not, discard $\varphi$. Else, $X=X_\varphi$, and  $(X_\varphi,\sigma)$ is a factor of $(X_\alpha,\sigma)$.
\end{procedure}

  The last step in this procedure is  algorithmic because of minimality and  Theorem 34 in \cite{CRS2012}. A computer program for this can be found at \cite{Walnut}. In some  cases the procedure can be executed by hand. We shall do this in Section \ref{sec: toeplitz} for the Toeplitz substitution, and in Section \ref{sec: morsefac} for the Thue-Morse substitution.

\medskip

 It is important to us that the last step in the procedure may be supplemented (and in many cases replaced) by checking whether  there exists an integer $p$ with $1\le p \le \mathrm{Card}(A_{\alpha}^{[3]})$ and an integer $M$ with $0\le M\le 2(L-1)$  , such that $\varphi^p$ is an amalgamation of
$(\widehat{\alpha}_{3,M})^p$, i.e., such that $\pi\circ (\widehat{\alpha}_{3,M})^p= \varphi^p\circ \pi$ holds for some letter-to-letter map $\pi$.

\smallskip

For an algorithm for the list problem for conjugacy we still need another ingredient. A dynamical system is called {\it coalescent} if
every endomorphism is an automorphism, i.e., every topological semi-conjugacy from the system onto itself is a  topological conjugacy.
 It was shown for a two symbol alphabet in \cite{Cov1971} and for a
general alphabet in \cite{Dur} that primitive, not necessarily constant length, substitutions generate coalescent dynamical systems.

\begin{procedure}\label{th: conj}
 Use Procedure~\ref{th: factor} to determine all primitive injective  substitutions $\beta$ with the same length that generate  factors of $(X_{\alpha},\sigma)$. Make the list for $\beta$, and check whether $\alpha$ is on it. If it is,  then $(X_{\alpha},\sigma)$ is  conjugate to $(X_{\beta},\sigma)$, by coalescence; if not, then  $(X_{\alpha},\sigma)$ is not conjugate to $(X_{\beta},\sigma)$. \qed
\end{procedure}

\section{The conjugacy class of the Toeplitz substitution}\label{sec: toeplitz}
We use Procedure~\ref{th: factor}  to determine the injective substitutions of length two that generate factors of the Toeplitz system $(X_\tau,\sigma)$ where $\tau$ is the substitution
$$ \tau(0)=01, \quad \tau(1)=00.$$
Actually, the property of $\tau$ that the first letters of the two $\tau$-blocks are equal implies  that for any $n$ $\tau^n(0)$ and $\tau^n(1)$ only differ in their final letter. It may be seen that it then suffices to restrict ourselves to 2-block  codes.

The set of  words of length two in ${\mathcal L}_\tau$ is equal to
$\mathcal{L}^2_\tau =\{00,01,10\},$ so we code the 2-blocks by  $A^{[2]}_\tau=\{1,2,3\}$.

We first consider the case where the letter-to-letter map $\pi$ is the identity. The graphs ${\mathcal G}_1={\mathcal G}_1^x,\; {\mathcal G}_{2,0}={\mathcal G}_{2,0}^x$ and ${\mathcal G}_{2,1}={\mathcal G}_{2,1}^x$ of a sequence $x$ in the 2-block presentation $X^{[2]}_\tau$ are given by

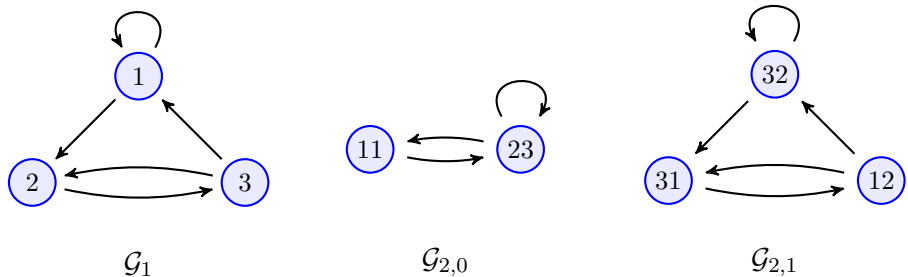
\begin{tikzpicture}[->,>=stealth',shorten <=8pt, shorten >=8pt,auto,node distance=2cm, thick,
  main node/.style={circle,fill=blue!10,draw,circle,inner sep=0pt},
  loop/.style={min distance=10mm,looseness=15,in=120,out=60},
  loop clock/.style={min distance=10mm,looseness=10,in=60,out=120}]
  \node[main node] (na) {$1$};
  \node[main node] [below left of=na] (nb) {$2$};
  \node[main node] [below right of=na] (nc) {$3$};
  \draw[blue,thick,fill=blue!8] (na) circle [radius=.31cm] (nb) circle [radius=.31cm] (nc) circle [radius=.31cm];
  \node at (na)  {$1$}; \node at (nb)  {$2$}; \node at (nc)  {$3$};
  \draw[->]
    (na) edge  (nb)
    (na) edge  [loop] node {} (na)
    (nb) edge [bend right=12] (nc)
    (nc) edge [bend right=12] (nb)
    (nc) edge  (na);
  \node  at (0,-2.5) {{\large ${\mathcal G}_1$}};
\end{tikzpicture}\label{G1-Toeplitz}\hspace*{.9cm}
\begin{tikzpicture}[->,>=stealth',shorten <=8pt, shorten >=8pt,auto,node distance=2cm, thick,
  main node/.style={circle,fill=blue!10,draw,circle,inner sep=0pt},
  loop/.style={min distance=10mm,looseness=10,in=120,out=60},
  loop clock/.style={min distance=10mm,looseness=10,in=60,out=120}]
  \node[main node] (na) {$11$};
  \node[main node] [right of=na] (nb) {$23$};
  \draw[blue,thick,fill=blue!8] (na) circle [radius=.31cm] (nb) circle [radius=.31cm];
  \node at (na)  {$11$}; \node at (nb)  {$23$};
  \draw[->]
    (na) edge [bend right=12] (nb) (nb) edge [bend right=12] (na)
    (nb) edge  [loop clock] node {} (nb);
  \node  at (1,-1.5) {{\large ${\mathcal G}_{2,0}$}};
\end{tikzpicture}\label{G20-Toeplitz}\hspace*{.9cm}
\begin{tikzpicture}[->,>=stealth',shorten <=8pt, shorten >=8pt,auto,node distance=2cm, thick,
  main node/.style={circle,fill=blue!10,draw,circle,inner sep=0pt},
  loop/.style={min distance=10mm,looseness=10,in=120,out=60},
  loop clock/.style={min distance=10mm,looseness=10,in=60,out=120}]
   \node[main node] (na) {$32$};
  \node[main node] [below left of=na] (nb) {$31$};
  \node[main node] [below right of=na] (nc) {$12$};
  \draw[blue,thick,fill=blue!8] (na) circle [radius=.31cm] (nb) circle [radius=.31cm] (nc) circle [radius=.31cm];
  \node at (na)  {$32$}; \node at (nb)  {$31$}; \node at (nc)  {$12$};
  \draw[->]
    (na) edge  (nb)
    (na) edge  [loop] node {} (na)
    (nb) edge [bend right=12] (nc)
    (nc) edge [bend right=12] (nb)
    (nc) edge  (na);
  \node  at (0,-2.5) {{\large ${\mathcal G}_{2,1}$}};
\end{tikzpicture}\label{G21-Toeplitz}

\medskip

There are two surjective graph homomorphisms $\varphi: \mathcal{G}_1\rightarrow\mathcal{G}_{2,0}$  which give a primitive substitution:
$$\varphi(1)=23,\; \varphi(2)=23,\; \varphi(3)=11,\quad \mathrm{and}\quad \varphi(1)=23,\; \varphi(2)=11,\; \varphi(3)=23.$$
The first $\varphi$ generates the 2-block presentation, since it may be checked that $\varphi=\widehat{\tau}_{2,0}$. After injectivization  it gives  the substitution $\alpha$ given by $\alpha(1)=13,\; \alpha(3)=11$, whose standard form is
the Toeplitz substitution. The second one is not equal to a $\widehat{\tau}_{2,M}$, and so we will postpone the answer to the question whether it generates a factor. It injectivizes to the substitution $\alpha$ given by $\alpha(1)=21,\; \alpha(2)=11$, which we call the \emph{rotated Toeplitz substitution}.

There is exactly one surjective graph homomorphism $\varphi: \mathcal{G}_1\rightarrow\mathcal{G}_{2,1}$,  which gives the primitive substitution:
\begin{equation*}\label{eq: LeMas}  \varphi(1)=32,\quad \varphi(2)=31,\quad \varphi(3)=12,\end{equation*}
which has the standard form given by $\alpha(1)=23,\; \alpha(2)=13, \; \alpha(3)=12.$ We call this substitution  {\it $3$-symbol Toeplitz}.
It may be checked that $\varphi=\widehat{\tau}_{2,1}$, and so the system generated by  this $\varphi$ is  conjugate to the Toeplitz system by Proposition \ref{th:thetahat}.

\smallskip

To finish, we  still have to examine the possibilities of letter-to-letter maps $\pi:\{1,2,3\}\rightarrow\{\dd,\ee\}$, where $\{\dd,\ee\}$ is a two letter alphabet.
There are three of these maps $\pi_k$ given by
$$\pi_1:\; 1\rightarrow\dd,\,2\rightarrow\dd,\,3\rightarrow\ee,\;\;\pi_2:\; 1\rightarrow\dd,\,2\rightarrow\ee,\,3\rightarrow\dd,\;\;
   \pi_3:\; 1\rightarrow\ee,\,2\rightarrow\dd,\,3\rightarrow\dd.$$
Let $t_k$ for $k=1,2,3$ be a sequence from $\pi_k(X_{\widehat{\tau}_{2,1}})$.
The graphs  ${\mathcal G}_1^1={\mathcal G}_1^{t_1},\; {\mathcal G}_{2,0}^1={\mathcal G}_{2,0}^{t_1}$ and ${\mathcal G}_{2,1}^1={\mathcal G}_{2,1}^{t_1}$
are given by

\smallskip

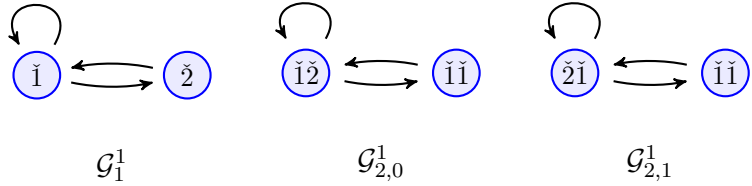
\begin{tikzpicture}[->,>=stealth',shorten <=8pt, shorten >=8pt,auto,node distance=2cm, thick,
  main node/.style={circle,fill=blue!10,draw,circle,inner sep=0pt},
  loop/.style={min distance=10mm,looseness=15,in=120,out=60},
  loop clock/.style={min distance=10mm,looseness=10,in=60,out=120}]
  \node[main node] (na) {$\dd$};
  \node[main node] [right of=na] (nb) {$\ee$};
  \draw[blue,thick,fill=blue!8] (na) circle [radius=.31cm] (nb) circle [radius=.31cm];
  \node at (na)  {$\dd$}; \node at (nb)  {$\ee$};
  \draw[->]
    (na) edge [bend right=12] (nb) (nb) edge [bend right=12] (na)
    (na) edge  [loop] node {} (na) ;
  \node  at (1,-1.2) {{\large ${\mathcal G}_1^1$}};
\end{tikzpicture}\label{G1-L2L-Toep}\hspace*{.5cm}
\begin{tikzpicture}[->,>=stealth',shorten <=8pt, shorten >=8pt,auto,node distance=2cm, thick,
  main node/.style={circle,fill=blue!10,draw,circle,inner sep=0pt},
  loop/.style={min distance=10mm,looseness=10,in=120,out=60},
  loop clock/.style={min distance=10mm,looseness=10,in=60,out=120}]
  \node[main node] (na) {$\dd\ee$};
  \node[main node] [right of=na] (nb) {$\dd\dd$};
  \draw[blue,thick,fill=blue!8] (na) circle [radius=.31cm] (nb) circle [radius=.31cm];
  \node at (na)  {$\dd\ee$}; \node at (nb)  {$\dd\dd$};
  \draw[->]
    (na) edge [bend right=12] (nb) (nb) edge [bend right=12] (na)
    (na) edge  [loop] node {} (na) ;
  \node  at (1,-1.2) {{\large ${\mathcal G}_{2,0}^1$}};
\end{tikzpicture}\label{G20-L2L-Toep}\hspace*{.5cm}
\begin{tikzpicture}[->,>=stealth',shorten <=8pt, shorten >=8pt,auto,node distance=2cm, thick,
  main node/.style={circle,fill=blue!10,draw,circle,inner sep=0pt},
  loop/.style={min distance=10mm,looseness=10,in=120,out=60},
  loop clock/.style={min distance=10mm,looseness=10,in=60,out=120}]
  \node[main node] (na) {$\ee\dd$};
  \node[main node] [right of=na] (nb) {$\dd\dd$};
  \draw[blue,thick,fill=blue!8] (na) circle [radius=.31cm] (nb) circle [radius=.31cm];
  \node at (na)  {$\ee\dd$}; \node at (nb)  {$\dd\dd$};
  \draw[->]
    (na) edge [bend right=12] (nb) (nb) edge [bend right=12] (na)
    (na) edge  [loop] node {} (na) ;
  \node  at (1,-1.2) {{\large ${\mathcal G}_{2,1}^1$}};
\end{tikzpicture}\label{G21-L2L-Toep}

There are obvious graph epimorphisms from ${\mathcal G}_1^1$ to ${\mathcal G}_{2,0}^1$ and to ${\mathcal G}_{2,1}^1$. The first one again yields
the Toeplitz substitution, the second one  yields the substitution
$$\check{\varphi}(\dd)=\ee\dd, \quad \check{\varphi}(\ee)=\dd\dd,$$
whose standard form is rotated Toeplitz. Since here we have the intertwining relation
$$\pi_1\circ\widehat{\tau}_{2,1}=\check{\varphi}\circ\pi_1,$$
$\check{\varphi}$ is an amalgamation of $\widehat{\tau}_{2,1}$, so $(X_{\check{\varphi}},\sigma)$ \emph{is} a  factor of the Toeplitz substitution system. It actually is conjugate to the Toeplitz system, since Toeplitz will be in the list of factors of the rotated Toeplitz substitution.

\smallskip

One can check that the
letter-to-letter map $\pi_2$ gives similar results, and that the graph ${\mathcal G}_1^3$ has two loops, which prevents graph homomorphisms in this case.

Conclusion: the conjugacy class of of the injective substitutions of the Toeplitz system consists of three substitutions:

\noindent Toeplitz, rotated Toeplitz, and 3-symbol Toeplitz.

\medskip

We will examine the properties of the minimal set $Y:=\pi_3(X_\tau^{[2]})\subseteq \{\dd,\ee\}$ in more detail. We showed that $Y$ is not generated by a substitution of length 2.
We will prove more: $Y$ is not generated by {\em any} substitution.
The only other example we know of this kind is the Rudin-Shapiro minimal set, cf.~\cite{Sch-Shallit}, page 1613.

First we prove the rather surprising fact that the sequences in $Y$ are essentially obtained by doubling the letters in the sequences of the Toeplitz minimal set.
Define the doubling morphism $\delta: \{0,1\}^*\rightarrow \{\dd,\ee\}^*$ by
$$\delta(0)=\dd\dd, \quad   \delta(1)=\ee\ee.$$

\begin{lemma} \label{lem:double} Let $\tau$ be the Toeplitz substitution on  $\{0,1\}$, let $\beta:= \widehat{\tau}_{2,1}$,   and let $\pi=\pi_3$ be the projection
$1\rightarrow\ee,\,2\rightarrow\dd,\,3\rightarrow\dd$.  Then   for all $n\ge 1$ \vspace*{-.15cm}
\begin{eqnarray*}\pi(\beta^{2n}(1))&=\ee\,\delta(\tau^{2n-1}(0))\,\ee^{-1},\\
  \pi(\beta^{2n}(2))&=\ee\,\delta(\tau^{2n-1}(1))\,\dd^{-1},\\
   \pi(\beta^{2n}(3))&=\dd\,\delta(\tau^{2n-1}(0)\,)\ee^{-1}.
\end{eqnarray*}
\end{lemma}
\noindent {\em Proof:}  By induction. For $n=1$ we have $\pi(\beta^{2}(1))=\pi(1231)=\ee\dd\dd\ee$. On the other hand,  $\ee\delta(\tau(0))\ee^{-1}=
\ee\delta(01)\ee^{-1}=\ee\dd\dd\ee\ee\ee^{-1}=\ee\dd\dd\ee$.  Now the induction step:
\begin{eqnarray*}
\pi(\beta^{2(n+1)}(1))&\!\!\!=\!\!\!&\pi(\beta^{2n}(1231))= \pi(\beta^{2n}(1))\pi(\beta^{2n}(2))\pi(\beta^{2n}(3))\pi(\beta^{2n}(1))\\   &\!\!\!=\!\!\!&\ee\delta(\tau^{2n-1}(0))\ee^{-1}\ee\delta(\tau^{2n-1}(1))\dd^{-1}\dd\delta(\tau^{2n-1}(0))\ee^{-1}\ee\delta(\tau^{2n-1}(0))\ee^{-1}\\
&\!\!\!=\!\!\!&\ee\,\delta(\tau^{2n-1}(0100))\ee^{-1}=\ee\,\delta(\tau^{2n+1}(0))\ee^{-1}.
\end{eqnarray*}
For the letters $2$ and $3$ a similar computation yields the corresponding formula \qed

\smallskip

It follows from Lemma \ref{lem:double} that $Y$ is the closed orbit of the sequence $y=\delta(t)$, where $t$ is the Toeplitz sequence.

 We need another combinatorial lemma.

\begin{lemma} Let $t$ be the Toeplitz sequence, and let $M$ be a fixed integer with $0\le M <2^n$ for some $n\ge 1$.
Then there is at most one word $w$ of length $2^n$ such that the square $ww$ occurs at some position $M\!\! \mod 2^n$ in $t$.
The same property holds for the sequence $y=\delta(t)$.\label{lem:square}
\end{lemma}

\noindent {\em Proof:} For even $n$ (for odd $n$  exchange the suffixes 0 and 1) the words
$$\tau^n(0)=:a_1a_2\dots a_{2^n\!-1}0,\, {\rm and}\; \tau^n(1)=:a_1a_2\dots a_{2^n\!-1}1$$
only differ in the last letter. Therefore the only two words of length $2^n$ occurring in $t$ at  position $M\!\! \mod 2^n$ are
$$ v_M:=a_{M+1}\dots a_{2^n\!-1}0a_1\dots a_M ,\, {\rm and}\; w_M:=a_{M+1}\dots a_{2^n\!-1}1a_1\dots a_M.$$
Since 11 does not occur in $t$, $\tau^n(11)$ does not occur in $t$, and this implies that $ v_M v_M$ is the only square occurring at  positions $M\!\! \mod 2^n$.

Now note that this implies that the same property holds for $\delta(t)$ for all words occurring at the even positions $2M\!\! \mod 2^n$.
But then it also holds for positions $2M\!+\!1\!\! \mod 2^n$, since if a square occurred at such an odd position, then we could shift 1 to the left,
obtaining a square at an even position (the words in $\delta(t)$ in even positions have prefix $\dd\dd$ or $\ee\ee$). \qed

We are now ready to prove the announced result.

\begin{proposition} Let $Y:=\pi_3(X_\tau^{[2]})\subseteq \{\dd,\ee\}^{\mathbb{Z}}$ be the projection of the 2-block presentation of the Toeplitz minimal set considered before. Then $(Y,\sigma)$ is not a substitution dynamical system.\label{prop:toeplitzdouble}
\end{proposition}

\noindent {\em Proof:} First note that if $Y$ would be generated by a substitution $\gamma$, then, by Cobham's Theorem, the length of $\gamma$ would be
a power of 2. Recall $y=\delta(t)$.  We use Lemma \ref{th: graph}. The graph  ${\mathcal G}_1^y$ is the complete graph on the nodes $\dd$ and $\ee$. For each $n$ and for all $M=0,...,2^n\!-1$ the graphs ${\mathcal G}_{2^n,M}^y$ have {\em only} one loop, because $y$ has only one square at position $M\!\! \mod 2^n$, by Lemma \ref{lem:square}. But then an epimorphism from ${\mathcal G}_1^y$ to ${\mathcal G}_{2^n,M}^y$ is impossible. \qed

\section{The length 2 substitution factors of the  Thue-Morse system}\label{sec: morsefac}
Let $\theta$ and $\theta^\flat$ be the Thue-Morse substitutions of length 2 on $A=\{0,1\}$ given by
$$ \theta(0)=01, \quad \theta(1)=10, \qquad \theta^\flat(0)=10,\quad \theta^\flat(1)=01.$$
The set of words of length 3 in the language of $\theta$ is
$\mathcal{L}_\theta^3=\{001,010,011,100,101,110\}.$
The usual lexicographic coding---which happens to be the binary coding---gives the 3-block alphabet $A_\theta^{[3]} :=\{1,2,3,4,5,6\}$.
The graph ${\mathcal G}_1={\mathcal G}_1^x$ of a sequence $x$ in the  3-block presentation $X_\theta^{[3]}$ is given by

\smallskip

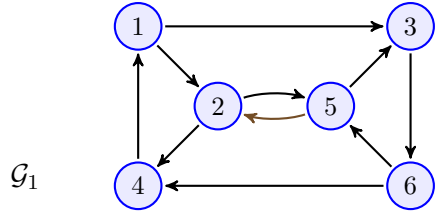
\hspace*{3cm}\begin{tikzpicture}[->,>=stealth',shorten <=8pt, shorten >=8pt,auto,node distance=3cm,
  thick,main node/.style={circle,fill=blue!10,draw,circle,inner sep=0pt,scale=0.5}]
  \node[main node] (n1) {$1$};
  \node[main node] [below right of=n1] (n2) {$2$};
  \node[main node] [right of=n2] (n5) {$5$};
  \node[main node] [above right of=n5] (n3) {$3$};
  \node[main node] [below left of=n2](n4) {$4$};
  \node[main node] [below right  of=n5](n6) {$6$};
  \draw[blue,thick,fill=blue!8]
   (n1) circle [radius=.31cm] (n2) circle [radius=.31cm] (n3) circle [radius=.31cm]
   (n4) circle [radius=.31cm] (n5) circle [radius=.31cm] (n6) circle [radius=.31cm];
  \node at (n1)  {$1$}; \node at (n2)  {$2$}; \node at (n3)  {$3$}; \node at (n4)  {$4$}; \node at (n5)  {$5$}; \node at (n6)  {$6$};
  \draw[->]
    (n4) edge (n1)
    (n2) edge (n4)
    (n1) edge (n2)
    (n1) edge (n3)
    (n6) edge (n4)
    (n3) edge (n6)
    (n5) edge (n3)
    (n6) edge (n5)
    (n5) edge [bend left=18,color=brown!60!black] (n2)
    (n2) edge [bend left=18] (n5);
    \node  at (-1.5,-2) {{\large ${\mathcal G}_1$}};
\end{tikzpicture}\label{G1-3-block-Morse}

\medskip

The graphs \;${\mathcal G}_{2,0}={\mathcal G}_{2,0}^x$ and ${\mathcal G}_{2,1}={\mathcal G}_{2,1}^x$ describing the 2-blocks in a sequence $x$ from the 3-block presentation $X^{[3]}_\theta$ are given by

\medskip

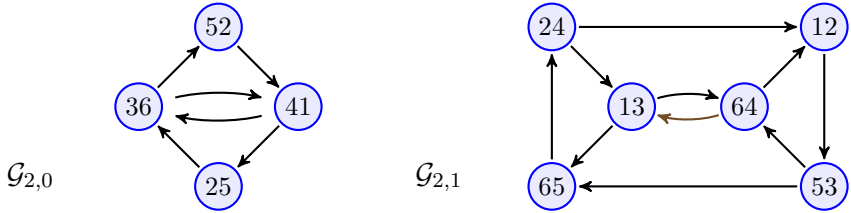
\begin{tikzpicture}[->,>=stealth',shorten <=8pt, shorten >=8pt,auto,node distance=3cm,
  thick,main node/.style={circle,fill=blue!10,draw,circle,inner sep=0pt,scale=0.5}]
  \node[main node] (n1) {$1$};
  \node[main node] [below left of=n1] (n2) {$2$};
  \node[main node] [below right of=n1] (n5) {$5$};
  \node[main node] [below left of=n5](n4) {$4$};
  \draw[blue,thick,fill=blue!8]
   (n1) circle [radius=.31cm] (n2) circle [radius=.31cm]
   (n4) circle [radius=.31cm] (n5) circle [radius=.31cm];
  \node at (n1)  {$52$}; \node at (n2)  {$36$}; \node at (n4)  {$25$}; \node at (n5)  {$41$};
  \draw[->]
    (n2) edge (n1)
    (n1) edge (n5)
    (n4) edge (n2)
    (n5) edge (n4)
    (n5) edge [bend left=14,shorten >=12pt,shorten <=12pt] (n2)
    (n2) edge [bend left=14,shorten >=12pt,shorten <=12pt] (n5);
    \node  at (-2.5,-2) {{\large ${\mathcal G}_{2,0}$}};
  \end{tikzpicture}\label{G20-3-block-Morse}\hspace*{1cm}
  \begin{tikzpicture}[->,>=stealth',shorten <=8pt, shorten >=8pt,auto,node distance=3cm,
  thick,main node/.style={circle,fill=blue!10,draw,circle,inner sep=0pt,scale=0.5}]
  \node[main node] (n1) {$1$};
  \node[main node] [below right of=n1] (n2) {$2$};
  \node[main node] [right of=n2] (n5) {$5$};
  \node[main node] [above right of=n5] (n3) {$3$};
  \node[main node] [below left of=n2](n4) {$4$};
  \node[main node] [below right  of=n5](n6) {$6$};
  \draw[blue,thick,fill=blue!8]
   (n1) circle [radius=.31cm] (n2) circle [radius=.31cm] (n3) circle [radius=.31cm]
   (n4) circle [radius=.31cm] (n5) circle [radius=.31cm] (n6) circle [radius=.31cm];
  \node at (n1)  {$24$}; \node at (n2)  {$13$}; \node at (n3)  {$12$}; \node at (n4)  {$65$}; \node at (n5)  {$64$}; \node at (n6)  {$53$};
  \draw[->]
    (n4) edge (n1)
    (n2) edge (n4)
    (n1) edge (n2)
    (n1) edge (n3)
    (n6) edge (n4)
    (n3) edge (n6)
    (n5) edge (n3)
    (n6) edge (n5)
    (n5) edge [bend left=18,color=brown!60!black] (n2)
    (n2) edge [bend left=18] (n5);
    \node  at (-1.5,-2) {{\large ${\mathcal G}_{2,1}$}};
\end{tikzpicture}\label{G21-3-block-Morse}

\medskip

To find all graph epimorphisms from ${\mathcal G}_1$ to ${\mathcal G}_{2,0}$ and ${\mathcal G}_{2,1}$, we exploit the following simple lemma.

\begin{lemma} Let $\varphi:{\mathcal G}\rightarrow {\mathcal G'}$ be a graph homomorphism. Suppose ${\mathcal G'}$ has no loops. Then 2-cycles and 3-cycles in  ${\mathcal G}$ are mapped to 2-cycles, respectively 3-cycles in ${\mathcal G'}$.\label{lem:loop}
\end{lemma}

It will appear that  all these graph epimorphisms are either a 3-block substitution of $\theta$ or of $\theta^\flat$.

We start with finding all $\varphi:{\mathcal G}_1\rightarrow {\mathcal G}_{2,0}$.
By  Lemma~\ref{lem:loop}, $\{\varphi(2),\varphi(5)\}$ equals $\{36,41\}$.
 If $\varphi(2)=36$, then $\varphi(4)=52$ and $\varphi(1)=41$, \emph{or}  $\varphi(4)=41$ and $\varphi(1)=25$.
 In the first case necessarily $(5,3,6) \rightarrow(41,25,36)$ by Lemma~\ref{lem:loop}, and we obtain
$$\widehat{\theta^\flat}_{\!2,2}: \quad  1\;\rightarrow41,\;2\rightarrow36,\;3\rightarrow25,\;4\rightarrow52,\;5\rightarrow41,\;6\rightarrow36.$$
In the second case $(5,3,6) \rightarrow(41,36,52)$, and we obtain
$$\widehat{\theta}_{2,0}: \quad  1\;\rightarrow25,\;2\rightarrow36,\;3\rightarrow36,\;4\rightarrow41,\;5\rightarrow41,\;6\rightarrow52.$$
 If $\varphi(2)=41$, then in the same way we obtain a third  and fourth epimorphism
 \begin{equation*}
  \begin{split}
 \widehat{\theta}_{2,2}:  \quad 1\rightarrow36,&\; 2\rightarrow41, 3\rightarrow52,\;4\rightarrow25,\;5\rightarrow36,\;6\rightarrow41, \\
 \widehat{\theta^\flat}_{\!2,0}: \quad 1\rightarrow52,&\; 2\rightarrow41, 3\rightarrow41,\;4\rightarrow36,\;5\rightarrow36,\;6\rightarrow25.
  \end{split}
\end{equation*}

Next we consider all $\varphi:{\mathcal G}_1\rightarrow {\mathcal G}_{2,1}$. Now $\{\varphi(2),\varphi(5)\}$  equals $\{13,64\}$.

\noindent If $\varphi(2)=13$, then $\varphi(4)=65$ and $\varphi(1)=24$, and also $\varphi(5)=64$, $\varphi(3)=12$ and $\varphi(6)=53$, since $(2,4,1)$ and $(5,3,6)$ form 3-cycles. In this way we obtain
$$\widehat{\theta^\flat}_{\!2,1}: \quad 1\rightarrow24,\; 2\rightarrow13, 3\rightarrow12,\;4\rightarrow65,\;5\rightarrow64,\;6\rightarrow53.$$
If $\varphi(2)=64$, then in the same way we obtain an epimorphism
$$\widehat{\theta}_{2,1}: \quad 1\;\rightarrow53,\;2\rightarrow64,\;3\rightarrow65,\;4\rightarrow12,\;5\rightarrow13,\;6\rightarrow24.$$

We now  do the letter-to-letter maps. This is much more involved than in the case of the Toeplitz substitution.

Note that the letter-to-letter maps from $A_\theta^{[3]}=\{1,2,3,4,5,6\}$ to another alphabet are in one to one correspondence with the set of all partitions
of $\{1,2,3,4,5,6\}$. Hence there are $B_6=203$ of such maps, where $B_6$ is the sixth Bernoulli number. Since $M$ can take the values 0 and 1, this means that there are 406 cases of candidate epimorphisms to consider.

To reduce this number, we note that there is the mirror symmetry $0\rightarrow1, 1\rightarrow 0$, which at the level of 3-blocks corresponds
to the permutation $\mathcal{P}= (16)(25)(34).$
Obviously a partition and its permuted version will  generate (if any) a substitution with the same standard form.

\noindent To further speed up the process we can apply the following three simple tools.

(T1) If ${\mathcal G}_{L,M}$ has more nodes than ${\mathcal G}_1$, then an epimorphism is not possible.

\noindent If there \emph{is} an epimorphism from ${\mathcal G}_1$ to ${\mathcal G}_{L,M}$, then:

(T2)  If the graph ${\mathcal G}_1$ contains a loop then ${\mathcal G}_{L,M}$ contains a loop.

(T3) If ${\mathcal G}_1$ and ${\mathcal G}_{L,M}$ have the same number of nodes, then they also must have the same  number of edges.


%

 With aid of the tools one finds 15 candidate substitutions to generate factors of the Thue-Morse system generated by injective substitutions of length 2:


\begin{center}
\hrule

{\footnotesize

\vspace*{.12cm}

\begin{tabular}{clcl}
{\large Nr.}\phantom{...}&  {\large  Partition}\phantom{...}&         {\large   $M$}\phantom{...} &                {\large Substitution }\\
\hline \\[-.3cm]
$\theta_1$\phantom{...}  & $ \{1,2,3\}\{4,5,6\}      $\phantom{...}& $0$\phantom{...} & $1\rightarrow 14,\:4\rightarrow 41          $ \\
$\theta_2$\phantom{...}  & $ \{1,2,3\}\{4,5,6\}      $\phantom{...}& $0$\phantom{...} & $1\rightarrow 41,\:4\rightarrow 14          $ \\
$\theta_3$\phantom{...} & $ \{1,2,5,6\}\{3,4\}       $\phantom{...}& $0$\phantom{...} & $1\rightarrow 31,\:3\rightarrow 11           $ \\
$\theta_4$\phantom{...}  & $ \{1,6\}\{2,3,4,5\}      $\phantom{...}& $0$\phantom{...} & $1\rightarrow 22,\:2\rightarrow 21          $ \\
$\theta_5$\phantom{...}  & $ \{1,4,5\}\{2,3\}\{6\}   $\phantom{...}& $1$\phantom{...} & $ 1\rightarrow 12 ,\:2\rightarrow 61,\: 6\rightarrow 21  $ \\
$\theta_6$\phantom{...}  & $ \{1,4,5\}\{2,6\}\{3\}   $\phantom{...}& $1$\phantom{...} & $1\rightarrow 21,\: 2\rightarrow 13,\: 3\rightarrow 12  $ \\
$\theta_7$\phantom{...}  & $ \{1,6\}\{2,5\}\{3,4\}   $\phantom{...}& $1$\phantom{...} & $1\rightarrow 23,\: 2\rightarrow 13,\: 3\rightarrow 12  $ \\
$\theta_8$\phantom{...}&$ \{1\}\{2,3\}\{4,5\}\{6\} $\phantom{...}& $0$\phantom{...} &$1\rightarrow 24,\:2\rightarrow 26,\: 4\rightarrow 41,\:6\rightarrow 42$ \\
$\theta_9$ \phantom{...}& $ \{1\}\{2,3\}\{4,5\}\{6\} $\phantom{...}& $0$\phantom{...} &$1\rightarrow 42,\:2\rightarrow 41,\: 4\rightarrow 26,\: 6\rightarrow 24$ \\
$\theta_{10}$\phantom{...}  & $ \{1,5\}\{2,6\}\{3\}\{4\}$\phantom{...}& $0$\phantom{...} &$1\rightarrow 41,\:2\rightarrow 32,\: 3\rightarrow 21,\: 4\rightarrow 12$\\
$\theta_{11}$\phantom{...}  & $ \{1,5\}\{2,6\}\{3\}\{4\}$\phantom{...}& $0$\phantom{...} &$1\rightarrow 32,\:2\rightarrow 41,\: 3\rightarrow 12,\: 4\rightarrow 21$\\
$\theta_{12}$\phantom{...}  & $ \{1,5\}\{2\}\{3\}\{4\}\{6\}   $\phantom{...}& $1$\phantom{...} &$1\rightarrow 13,\:2\rightarrow 64,\: 3\rightarrow 61,\: 4\rightarrow 12,\: 6\rightarrow 24 $ \\
$\theta_{13}$\phantom{...}& $ \{1\}\{2,3\}\{4\}\{5\}\{6\}$\phantom{...}& $1$\phantom{...} &$1\rightarrow 24,\:2\rightarrow 12,\: 4\rightarrow 65,\: 5\rightarrow 64,\: 6\rightarrow 52 $ \\
$\theta_{14}$\phantom{...}  & $  \{1\}\{2\}\{3\}\{4\}\{5\}\{6\} $\phantom{...}& $1$\phantom{...} &$1\rightarrow 24,\: 2\rightarrow 13,\: 3\rightarrow 12,\:4\rightarrow 65,\:5\rightarrow 64,\:6\rightarrow 53$\\
$\theta_{15}$\phantom{...}  & $ \{1\}\{2\}\{3\}\{4\}\{5\}\{6\}  $\phantom{...}& $1$\phantom{...} & $1\rightarrow 53,\: 2\rightarrow 64,\: 3\rightarrow 65,\: 4\rightarrow 12,\: 5\rightarrow 13,\: 6\rightarrow 24$
\end{tabular}}

\smallskip

{ Thue-Morse Factor List---direct projections}

\smallskip

\hrule

\smallskip

\end{center}

All 15 do generate a factor by the following arguments. 
The systems generated by  $\theta_1, \theta_2, \theta_3$ and $\theta_4$ are well known factors of the Thue-Morse system, and  $\theta_8$, $\theta_9, \theta_{10}, \theta_{11}, \theta_{14}$, and $\theta_{15}$  actually give conjugate systems, because they are injectivizations of 3-block substitutions of $\theta$ or of $\theta^\flat$.
 All others turn out to be amalgamations of either $\theta_{14}$ or $\theta_{15}$. For example $\theta_5\circ  \pi =\pi\circ \theta_{15}$, where the partition representation of $\pi$ is $\{1,4,5\}\{2,3\}\{6\}$. In the same way $\theta_6, \theta_7, \theta_{12}$ and $\theta_{13}$  are amalgamations of respectively
$\theta_{14}, \theta_{14}, \theta_{15}$ and $\theta_{14}$ by projections whose partition representation can be found in the table.

\smallskip

At an early stage of our research we found more than 15 substitutions in the factor list, failing to see that some were essentially the same.
For example, $\theta_8$ can also be  obtained   as the substitution generated by  the partition $\{1,5\}\{2,6\}\{3\}\{4\}$, but now for $M=1$.
It is therefore important to transform all $\theta_{k}$ to their standard forms $\thetaSF_k$.
The  standard forms of the substitutions in the Thue-Morse factor list are  given in the following table.

\smallskip

\begin{center}
\hrule

{\footnotesize

\vspace*{.1cm}

\begin{tabular}{clcl}
{\large Nr.} & {\large Standard form} & \hspace{-1cm} {\large Nr.} & {\large Standard form}\\
\hline \\[-.3cm]
$\thetaSF_1$& $1\rightarrow 12,\:2\rightarrow 21 $ & \hspace{-1cm} $\thetaSF_9 \;$&$1\rightarrow 23,\:2\rightarrow 14,\: 3\rightarrow 21,\: 4\rightarrow 12$ \\
$\thetaSF_2$& $1\rightarrow 21,\:2\rightarrow 12 $ & \hspace{-1cm} $\thetaSF_{10}$&$1\rightarrow 21,\:2\rightarrow 13,\: 3\rightarrow 43,\:4\rightarrow 31$ \\
$\thetaSF_3$& $1\rightarrow 21,\:2\rightarrow 11 $ & \hspace{-1cm} $\thetaSF_{11}$&$1\rightarrow 23,\:2\rightarrow 13,\: 3\rightarrow 41,\: 4\rightarrow 31$\\
$\thetaSF_4$& $1\rightarrow 12,\:2\rightarrow 11 $ & \hspace{-1cm} $\thetaSF_{12}$&$1\rightarrow 12,\:2\rightarrow 31,\: 3\rightarrow 45,\: 4\rightarrow 35,\: 5\rightarrow 14 $ \\
$\thetaSF_5$& $ 1\rightarrow 12,
\:2\rightarrow 31,\: 3\rightarrow 21  $ & \hspace{-1cm} $\thetaSF_{13}$&$1\rightarrow 21,\:2\rightarrow 13,\: 3\rightarrow 45,\:  4\rightarrow 51,\: 5\rightarrow 43 $ \\
$\thetaSF_6$& $1\rightarrow 21,
\: 2\rightarrow 13,\: 3\rightarrow 12  $ &\hspace{-1cm} $\thetaSF_{14}$&$1\rightarrow 23,\: 2\rightarrow 14,\: 3\rightarrow  21,\:4\rightarrow 56,\:5\rightarrow 63,\:6\rightarrow 54$\\
$\thetaSF_7$& $1\rightarrow 23,
\: 2\rightarrow 13,\: 3\rightarrow 12  $ &\hspace{-1cm} $\thetaSF_{15}$&$1\rightarrow 23,\: 2\rightarrow 13,\: 3\rightarrow 41,\: 4\rightarrow 56,\: 5\rightarrow 46,\: 6\rightarrow 25$\\
$\thetaSF_8$& $1\rightarrow 12,\:2\rightarrow 31,\: 3\rightarrow  34,\:4\rightarrow 13$ & &
\end{tabular}}

\hrule

\end{center}

\bigskip

\section{The Thue-Morse conjugacy list}\label{sec: morse}

Three substitutions  ($\theta_3$, $\theta_4$ and $\theta_7$) in the Thue-Morse factor list  generate systems that are certainly not conjugate to the Thue-Morse system, as they are in the Toeplitz conjugacy class. Obviously $\theta_1$ and  $\theta_2$ are in the conjugacy list, and we already know that the substitutions $\theta_8$, $\theta_9, \theta_{10}, \theta_{11}, \theta_{14}$, and $\theta_{15}$  generate  systems conjugate to the Thue-Morse system.
To see whether the 4 remaining substitutions yield systems  conjugate to the Thue-Morse system, according to Procedure~\ref{th: conj} we would have to construct the factor list of each of these. This is quite involved, for example the 3-block presentations of the two factors on 5 symbols have 11 symbols.

However, there is a quicker way to determine whether these factors are conjugate to the Thue-Morse system, by finding explicit semi-conjugacies from these factors to the Thue-Morse system. Then by coalescence the systems are conjugate.

It is quickly verified that indeed each of  $\theta_5, \theta_6, \theta_{12}$ and $\theta_{13}$  amalgamates to Morse or Morse flat. For example for $\theta_{12}$ one takes $1,4\rightarrow 0,\; 2,3,6\rightarrow 1.$

Conclusion: there are 12 primitive injective substitutions of length 2 that generate a system conjugate to the Thue-Morse dynamical system.

\section{Proper factors}\label{sec: proper}

We have seen that for the Thue-Morse system all factors are actually conjugate to the system, if there are no spectral obstructions. In this section we present in a simple way a system with mixed spectrum which has another system with mixed spectrum as a proper factor.

Let $\alpha$ be the Mephisto Waltz substitution given by $\alpha(1)=112,\quad \alpha(2)=221.$
Let $\beta$ be the substitution on four symbols given by $$\beta(1)=123,\quad \beta(2)=124,\quad\beta(3)=341,\quad \beta(4)=431.$$

\begin{proposition}\label{prop: devil} The system $(X_\alpha,\sigma)$ is a proper factor of $(X_\beta,\sigma)$.
\end{proposition}

 \noindent {\it Proof:}
Note that $\alpha$ is an amalgamation of $\beta$ under the projection map $$\pi(1)=\pi(2)=1,\quad \pi(3)=\pi(4)=2.$$ Therefore $(X_\alpha,\sigma)$ is a factor of $(X_\beta,\sigma)$. However, $(X_\beta,\sigma)$ is not a factor of $(X_\alpha,\sigma)$. To see this, note that 13 and 14 are
 in $\mathcal{L}_\beta$, and that 1 is suffix of $\beta^2(1)$. It follows that the two sequences $z:=(\beta^2)^{\infty}(1)\cdot(\beta^2)^{\infty}(3)$ and $z':=(\beta^2)^{\infty}(1)\cdot(\beta^2)^{\infty}(4)$ are in  $X_\beta$.
  Next, note that $z\ne z'$ and that $\pi(z)=\pi(z')$, since for all $n$
$$\pi \beta^n(3)=\alpha^n\pi(3)=\alpha^n(2)=\alpha^n\pi(4)=\pi \beta^n(4).$$
 Now suppose $\psi:X_\alpha\rightarrow X_\beta$ is a semi-conjugacy. Then, by coalescence, $\pi\circ\psi$ is a conjugacy. But this contradicts our finding that $\pi$ is 2 to 1 somewhere.\qed

 \smallskip

We remark that it is quite a delicate matter whether a factor is proper or not. For example, let $\alpha$ be the Mephisto Waltz, and  let $\delta$ be the substitution  defined by
$$\delta(1)=123,\quad \delta(2)=124,\quad\delta(3)=431,\quad \delta(4)=432. $$
Then $(X_\delta,\sigma)$ is conjugate to $(X_\alpha,\sigma)$, since it may be easily checked that  $\delta$ is the injectivization of the 3-block substitution $\widehat{\alpha}_{3,0}$.

 However, suppose we would follow the approach above, noting that $\alpha$ is an amalgamation of $\delta$ with the same $\pi$ map as above. Now $\delta^2$ has fixed prefix letters 1 and 4 and fixed suffix letters 1,2,3 and 4. This implies that the eight sequences $z_{b,a}:=(\delta^2)^{\infty}(b)\cdot(\delta^2)^{\infty}(a)$ are well-defined for $a=1,4$ and $b=1,2,3,4$. But, similarly as above, we have $\pi(z_{b,1})=\pi(z_{b',1})$ and $z_{b,1}\ne z_{b',1}$ for $b=1, b'=2$ and for $b=3, b'=4$ , yielding several points where $\pi$ is 2 to 1. However, this does not contradict conjugacy of the two systems, since neither $z_{2,1}$ nor $z_{4,1}$ are  elements of $X_\delta$, simply because the words 21 and 41 are not in the language of $\delta$.

 \section{Epilogue}\label{sec: epilogue}

   Related work can be found in the thesis of Joseph Herning  \cite{Hern} which mainly concentrates on bijective substitutions, which generate a relatively small subclass of systems with partially continuous spectrum. A {\em bijective substitution} $\alpha$ on an alphabet $A$ is defined by $\{\alpha(a)_i: a\in A\}=A$ for all $1\le i\le L$. One of the major results in \cite{Hern} is that there exist substitution dynamical systems that do not have discrete spectrum factors generated by substitutions. As an example Herning gives the substitution $\alpha$ on three symbols, which also occurs in \cite{Que}, defined by
   $$\alpha(1)=121, \; \alpha(2)=233, \; \alpha(3)=312.$$
   We have reproved his result by computing the factor list of $\alpha$. It consists of nine injective substitutions, on alphabets of size three to eight, all (indeed!) generating systems with partially continuous spectrum. Without doing any computations, it follows from Theorem 8 in \cite{Mentz} that these factors are in fact all conjugate to the system generated by $\alpha$, since the substitution $\alpha$ has no non-trivial amalgamations.

\smallskip

An interesting extension of our result would be to consider also non-constant length substitutions. For example, let $\theta$ be the ternary Thue-Morse substitution, defined by
$$\theta(1)=123,\; \theta(2)=13,\; \theta(3)=2.$$
An application of Theorem 1 in Section V of \cite{Dek78} shows that $(X_\theta,\sigma)$ is conjugate to a substitution of constant length 2 on 6 symbols. Its injectivization  is a substitution on 5 symbols, and taking the standard form of this substitution we find that it is on the Thue-Morse list.

\smallskip

The paper \cite{Salo} considers conjugacies between systems generated by two primitive substitutions whose matrices have the same Perron-Frobenius eigenvalue: it is shown there that modulo powers of the shift there are only finitely many conjugacies between such systems. Nevertheless, it has  been shown in \cite{Dek2014} that there are infinitely many systems on the Thue-Morse list, all generated by primitive injective substitutions with Perron-Frobenius eigenvalue 2.

\smallskip

Primitive substitutions generate dynamical systems with a unique shift invariant measure. One can consider Problem~\ref{prob:b} for measure-theoretic conjugacy.
When a substitution of length $L$ generates a system with discrete spectrum, then obviously there are infinitely many primitive injective substitutions in the measure-theoretic conjugacy class (in fact {\it all} pure (see \cite{Dek78})  substitutions of length $L$) . When there is partially continuous spectrum, we believe that the equivalence class will be finite, and the same as for topological conjugacy. This has been proved for a subclass of such constant length substitutions in \cite{HP}.

\section{Acknowledgement}
We have profited from electronic discussions  on coalescence with Fabien Durand and Reem Yassawi, and from useful remarks by Michelle Lemasurier.

 \bibliographystyle{plain}

 \end{document}